\documentclass{article}
\usepackage{url}
\usepackage{array}
\newtheorem{lemma}{Lemma}

\newtheorem{theorem}{Theorem}
\newtheorem{conjecture}{Conjecture}

\newtheorem{preproof}{{\bf Proof.}}

\newenvironment{proof}[1]{\begin{preproof}{\rm
               #1}\hfill{\rule[-0.5mm]{2mm}{2mm}}}{\end{preproof}}

\usepackage{array}
\newlength{\cellwid}
\makeatletter
\newenvironment{latinsq}[1][00]{%
 \def\centcol{\centering \let\\=\tabularnewline \CellStrut}
 \arraycolsep0pt
 \setbox\@tempboxa\hbox{#1}%
 \cellwid\ht\@tempboxa  \advance\cellwid\dp\strutbox
 \ifdim\cellwid<\wd\@tempboxa 
   \@tempdima .5\wd\@tempboxa \advance\@tempdima -.5\cellwid 
   \advance\cellwid \@tempdima \advance\@tempdima\dp\strutbox
 \else
   \@tempdima\dp\strutbox
 \fi
 \edef\CellStrut{\vrule
   width\z@ height\the\cellwid depth\the\@tempdima \relax}
 \advance\cellwid\@tempdima \advance\cellwid-2\arraycolsep
 \array{|>{\centcol}p{\cellwid}|*{20}{>{\centcol}p{\cellwid}|}}}%
 {\endarray}
\makeatother

\title{Critical sets in the elementary abelian 2- and 3- groups}
\author{Richard Bean}
\date{}

\begin{document}
\maketitle

\begin{center}
Institute for Studies in Theoretical Physics and Mathematics \\
PO Box 19395-5746 \\
Tehran, I.R. Iran \\
\end{center}

\begin{abstract}

In 1998, Khodkar showed that the minimal critical set in the Latin
square corresponding to the elementary abelian 2-group of order 16
is of size at most 124.  Since the paper was published, improved
methods for solving integer programming problems have been
developed.  Here we give an example of a critical set of size 121
in this Latin square, found through such methods. We also give a
new upper bound on the size of critical sets of minimal size for
the elementary abelian 2-group of order $2^n$:
$4^{n}-3^{n}+4-2^{n}-2^{n-2}$. We speculate about possible lower
bounds for this value, given some other results for the elementary
abelian 2-groups of orders 32 and 64. An example of a critical set
of size 29 in the Latin square corresponding to the elementary
abelian 3-group of order 9 is given, and it is shown that any such
critical set must be of size at least 24, improving the bound of
21 given by Donovan, Cooper, Nott and Seberry.

\end{abstract}

\section{Definitions}
A Latin square $L$ of order $n$ is an $n \times n$ array of
entries $\{ (i,j;k) \}$ such that each row and column of $L$
contains each of $n$ possible elements exactly once. In this paper
these $n$ possible elements will be the set $\{ 0, \dots, n-1 \}$,
as we will be using bit-wise logical operators later.  A subsquare
$S$ of order $m$ from $L$ is an array of $m \times m$ entries from
$L$, not necessarily contiguous, such that each row and column of
$S$ contains each of the $m$ possible elements exactly once.

A {\it critical set} $C$ of $L$ is a subset of $L$ such that $L$
is the only superset of $C$ which is a Latin square (is {\it
uniquely completable}) and no subset of $C$ has this property.

We now review the definition of a Latin trade (see for example
\cite{closing}). Let $P$ and $P'$ be two partial Latin squares of
the same order, with the same size and shape. Then $P$ are $P'$
are said to be {\it mutually balanced} if the entries in each row
(and column) of $P$ are the same as those in the corresponding row
(and column) of $P'$. They are said to be {\it disjoint}
 if no cell in $P'$ contains the same entry as the corresponding cell
in $P$. A {\it Latin trade} (also known as a {\it Latin
interchange}) $I$ is a partial Latin square for which there exists
another partial Latin square $I'$, of the same order, size and
shape with the property that $I$ and $I'$ are disjoint and
mutually balanced. The partial Latin square $I'$ is said to be a
{\it disjoint mate} of $I$.  A Latin trade of size 4 is known as
an {\it intercalate}. An example of a Latin trade is given below.

$$
\begin{tabular}{cc}
$\begin{array}{|c|c|c|c|c|c|}
\hline 9 & . & . & . & . & 3 \\
\hline 3 & 4 & 5 & 6 & 7 & 8 \\
\hline 4 & 5 & 6 & 7 & 8 & 9 \\
\hline \end{array}$ \quad\quad
$\begin{array}{|c|c|c|c|c|c|}
\hline 3 & . & . & . & . & 9 \\
\hline 4 & 5 & 6 & 7 & 8 & 3 \\
\hline 9 & 4 & 5 & 6 & 7 & 8 \\
\hline \end{array}$
\end{tabular}
$$
\hspace{44mm} $I$ \hspace{38mm} $I'$ \newline

The following lemma clarifies the connection between critical sets and Latin
trades.

\begin{lemma}
A partial Latin square $C\subset L$, of size $s$ and
order $n$, is a critical set for a Latin square $L$ if and only if the
following hold:
\begin{enumerate}
\item $C$ contains an element of every Latin trade that occurs in $L$;
\item for each $(i,\, j;\, k)\in C$, there exists a Latin trade
$I_r$ in $L$ such that $I_r\cap C = \{(i,\, j;\, k)\}.$
\end{enumerate}
\end{lemma}

\begin{proof}{ \noindent
\begin{enumerate}
\item If $C$ does not contain an element from some Latin trade
$I$, where $I$ has  disjoint mate $I'$, then $C$ is also a partial
Latin square of $L'= (L\setminus I)\cup I'$. Hence $C$ is not uniquely
completable.
\item If no such Latin trade $I_r$ can be found, then the position
$(i,\, j;\, k)$ may be deleted from $C$ and $C\setminus\{(i,\,j;\,k)\}$
will still be uniquely
completable and thus a critical set for $L$.
\end{enumerate}
}\end{proof}

The elementary abelian m-group of order $m^n$ is the group $Z_m$
multiplied with itself $n$ times, denoted in this paper by
$L(m^{n})$.  Elsewhere, it has also been denoted by $Z_m^n$.

\section{Integer Programming}
To find a critical set $C$ in the elementary abelian 2-group of
order $2^{n}$, it thus suffices to solve an integer programming
problem.  If $\mathcal{T}$ is the complete set of Latin trades in
a Latin square $L$, then the optimal solution to the following
integer program is the size of the minimal critical set in $L$.

\begin{eqnarray*}
{\rm Minimize}: & \sum_{x \in L}{C_{x}} & \\
{\rm Subject~to}: \\
{\rm For~each} &T \in \mathcal{T},& \sum_{x \in T}{C_{x}} \geq 1 \\
\end{eqnarray*}
Throughout this paper, $C_{x}$ is 1 if $x \in C$ and 0 otherwise.

In \cite{MR99f:05019} Khodkar uses this method to prove that the
minimal critical set for $L(8)$ is 25 and finds a critical set of
size 124 in $L(16)$.  Even for $L(8)$, it is a difficult task to
find all the Latin trades in this square, and so a subset of the
Latin trades must be used.  For $L(8)$ Khodkar chose all the
subsquares isomorphic to $L(2)$ and $L(4)$ because a linear
program similar to the above can be written as follows. $I_1$
denotes all the subsquares isomorphic to $L(2)$ in $L(8)$ and
$I_2$ denotes all the subsquares isomorphic to $L(4)$ in $L(8)$.

\begin{eqnarray*}
{\rm Minimize}:  & \sum_{x \in L(8)}{C_{x}} & \\
{\rm Subject~to}:  \\
{\rm For~each} &T \in {I_{1}},& \sum_{x \in T}{C_{x}} \geq 1 \\
{\rm For~each} &T \in {I_{2}},& \sum_{x \in T}{C_{x}} \geq 5 \\
\end{eqnarray*}

As explained by Khodkar, although this does not cover all the
Latin trades in $L(8)$, this works because the inequalities
above make each subsquare isomorphic to $L(4)$ in $L(8)$
uniquely completable.  Since a critical set of size 25 can be found
in $L(8)$ and the optimal solution to the above integer programming
problem is also of size 25, we know that the minimal critical set
size in $L(8)$ is also 25.

Khodkar used the CPLEX solver \cite{cplex}, version 6.5, for solving these
integer programming problems.  For solving integer programs,
CPLEX uses a branch-and-bound or branch-and-cut approach.  Although
CPLEX can find provably optimal solutions, for finding ``good'' solutions
quickly a non-deterministic search is much better.  We use this approach
to find a critical set in $L(16)$ and much better solutions to
similar integer programs for $L(32)$ and $L(64)$.

Walser extended the WalkSAT algorithm, a satisfiability (SAT)
solver using local search \cite{sel}, from solving SAT problems to
solving integer linear programming problems.  The resulting
program is called WSAT(OIP) \cite{wsat}.  We consider the
following problem, also from the Khodkar paper, where $I_n$ is the
set of all subsquares in $L(16)$ isomorphic to $L(2^n)$.
\newline

\begin{eqnarray*}
{\rm Minimize}:  & \sum_{x \in L(16)}{C_{x}} & \\
{\rm Subject~to}:  \\
{\rm For~each} &T \in {I_{1}},& \sum_{x \in T}{C_{x}} \geq 1 \\
{\rm For~each} &T \in {I_{2}},& \sum_{x \in T}{C_{x}} \geq 5 \\
{\rm For~each} &T \in {I_{3}},& \sum_{x \in T}{C_{x}} \geq 25 \\
\end{eqnarray*}

Khodkar found a solution of size 124 to the above problem with
CPLEX which was a critical set for $L(16)$.  He noted \cite{ak}
that CPLEX did find smaller solutions which were not critical
sets.  This is the case because the inequalities above do not
necessarily force the subsquares isomorphic to $L(8)$ to be
uniquely completable, though they force subsquares isomorphic to
$L(4)$ and $L(2)$ to be uniquely completable.  Also, other Latin
trades not specified in the program are missing from $L(16)$.

Using the CPLEX 8.0 solver on this problem, it took several hours
on an Pentium-4 system at 1.5 Ghz, with a variety of options, to
find a solution of size 118.  In contrast WSAT 1.105 can find
solutions of size 112 to the above problem in a few seconds.  The
disadvantage is that it cannot prove such solutions are optimal.
We note the right hand side (RHS) of the inequalities for all such
solutions of size 112 found are as follows: 512 with RHS 1, 208
with RHS 2 or 3, 32 with RHS 4; 72 with RHS 5 or 7, 40 with RHS 9,
8 with RHS 10, 11 or 12; 28 with RHS 28; and 16 with RHS 26 or 30.
Such symmetry and universal divisibility by powers of 2 leads the
author to conjecture that 112 is in fact the optimal solution to
the above problem.  Also, in all such solutions of size 112, each
of the subsquares isomorphic to $L(8)$ is uniquely completable.

It is simple to find critical sets of size 121 using this method.
Here is an example:

\begin{center}
$
\begin{latinsq}[2]
\hline 0 && 2 &&&& 6 & 7 & 8 & 9 &&& 12 &&& 15   \\
\hline & 0 &&& 5 & 4 & 7 && 9 & 8 && 10 & 13 & 12 &&  \\
\hline 2 & 3 & 0 &&&& 4 &&&&&&&&& 13   \\
\hline 3 & 2 &&& 7 &&&& 11 &&& 8 && 14 && 12   \\
\hline && 6 &&& 1 &&& 12 & 13 && 15 & 8 && 10 & 11   \\
\hline 5 && 7 && 1 & 0 && 2 && 12 & 15 & 14 &&&&  \\
\hline 6 & 7 & 4 & 5 & 2 &&&&&& 12 && 10 & 11 &&  \\
\hline && 5 & 4 &&&& 0 & 15 &&& 12 & 11 & 10 & 9 &  \\
\hline 8 &&& 11 && 13 & 14 &&& 1 &&& 4 & 5 & 6 &  \\
\hline && 11 &&& 12 &&& 1 & 0 & 3 && 5 & 4 && 6   \\
\hline &&& 9 & 14 & 15 &&& 2 &&& 1 &&& 4 &  \\
\hline & 10 &&&& 14 & 13 & 12 &&& 1 & 0 && 6 &&  \\
\hline 12 & 13 && 15 &&&&&& 5 & 6 & 7 & 0 &&& 3   \\
\hline && 15 & 14 & 9 && 11 & 10 &&& 7 &&&& 3 & 2  \\
\hline 14 &&& 13 & 10 && 8 & 9 && 7 && 5 &&&&  \\
\hline &&&& 11 && 9 & 8 & 7 & 6 & 5 && 3 && 1 &  \\
\hline
\end{latinsq}
$
\end{center}

Methods to try to improve this bound such as modifying the RHS
bounds of 5 and 25 to various other values ranging from 6 or 7 for
the $I_2$ trades and 26 to 30 for the $I_3$ trades have been
unsuccessful. Similarly finding the missing trades in solutions of
size 25 or more to the third set of inequalities has not aided
much.

If constraints are added to the program to ensure that each
row and column contains either 7 or 8 entries, and each element
occurs either 7 or 8 times, a critical set of size 121 can
still be found.  However this is still the best known solution.
To improve this result, the author suggests that a number of methods
be tried.

\begin{itemize}
\item The structure of the problem best suits a 64-bit CPU - any
subset of $L(16)$ can be represented as 4 64-bit words.  To reduce
the size of a known critical set, perhaps the same method as
\cite{closing} could be used - for increasing $n$, attempt to
remove $n+1$ entries while adding $n$ entries, and check for
unique completion afterwards.

\item An approach similar to the paper \cite{sts} which concerned
Steiner triple systems can be used.  That is, if we find that our
integer program solver outputs a solution which is not uniquely
completable, we determine which trades in $L(16)$ are ``missing''
from the solution, use the automorphism group of $L(16)$ to
generate copies of all such trades, ``minimize'' this list, add
them to the integer program, and continue.

\item Recently, Margot \cite{margot} wrote about pruning 0-1
integer programs which have a high degree of symmetry, and used
this method to prove that no 4-(10,5,1) covering design with less
than 51 sets exists.  Although the problem for $L(16)$ may have
too many variables for an optimal solution to be found, there is
hope that this method may find a better solution.  Similarly, the
integer programming methods of Applegate, Raines, and Sloane
\cite{apple} may be applicable here.
\end{itemize}

Applying the above integer program to $L(32)$ the best solution
found so far has been of size 546.  Using just intercalates, a
solution of size 514 has been found.  It has not been
computationally feasible to determine whether or not either
solution has unique completion.  The best bound currently known is
obtained by a method similar to that of Donovan, Fu and Khodkar
\cite{dfk}; we take the critical set of size $4^5-3^5=781$ as
constructed by Stinson and van Rees \cite{svr}, add in the main
backwards diagonal and the last row and column, and recursively
remove entries beginning in the top left-hand corner.  This
results in a critical set of size 658.

Given the evidence above, $L(32)$ seems not to have a critical set
of size less than or equal to $\frac{32^2}{2}$.  Ghandehari,
Hatami, and Mahmoodian \cite{ghm} found that there exists a Latin
square $L$ such that the minimal critical set in $L$ has size
$n^2-(e+o(1))n^{5/3}$.  The author conjectures that the smallest
order $n$ for which there exists a Latin square $L$ of order $n$
with no critical set of size less than $\frac{n^2}{2}$ is $n=32$,
in the square $L(32)$.

Similarly, the best solution to the above integer program for
$L(64)$ has been of size 2470.

With these results in mind, we conjecture the following bound.

\begin{conjecture}
If $C$ is a critical set of minimal size size in $L(2^n)$, $n \geq
4$, $121.4^{n-4} \leq |S| \leq 5^{n-1}$.
\end{conjecture}

\section{A new bound for the size of a minimal critical set in $L(2^n)$}
We can find a new general bound on the size of critical sets in
$L(2^n)$ by taking the same Stinson and van Rees construction for
a critical set in $L(2^n)$, adding 3 entries and deleting
$2^n+2^{n-2}-1$ entries.  This results in another critical set.

\begin{theorem}
If $C$ is a critical set of minimal size in $L(2^n)$, $|C| \leq
4^n-3^n+4-2^n-2^{n-2}$.
\end{theorem}

\begin{proof}{
We describe a construction isomorphic to the
one given by Stinson and van Rees for a critical set $C$ in
$L(2^n)$. Formally, if we consider the critical set $C$ in
$L(2^n)$, the entry at $(i,j) \in C$ if and only if ($i$ \& $j)
\neq 0$, where $\&$ is the bitwise logical AND operator.  Since we
are dealing with $L(2^n)$, the element at $(i,j)$ is $i \wedge j$
where $\wedge$ is the logical XOR operator.

Thus $C = \{ (i,j;i \wedge j) | (i\&j) \neq 0 \}$.  Consider $D=
(C \cup \{ (0,0;0),$ $(0,2^n-1;2^n-1),$ $(2^n-1,0;2^n-1) \})
\setminus (\{ (2^n-1, x; 2^n - 1 - x), (x, 2^n-1; 2^n - 1 - x) |
2^{n-1} \leq x \leq 2^n-1 \} \cup \{ (x, x; 0) | 2^{n-2} \leq x
\leq 2^{n-1}-1 \}$).

$D$ has unique completion as we first add back in the entries $\{
(x,x;0) | 2^{n-2} \leq x \leq 2^{n-1}-1 \}$, followed by the entry
$(2^n-1, 2^n-1; 0)$, then the entries $\{ (2^n-1, x; 2^n-1-x), (x,
2^n-1; 2^n-1-x) | 2^{n-1} \leq x < 2^n-1 \}$ in any order,
resulting in a superset of $C$.

To prove that every entry in $C$ is necessary, we say that an
intercalate $I \subseteq L(2^n)$ {\it proves the necessity} of an
entry $x \in C$ if $I \cap C = \{x\}$.  If we can find such trades
for every entry in $C$, then $C$ is a critical set, by Lemma 1.

We prove the necessity of the two additional entries
$(0,2^n-1;2^n-1), (2^n-1,0;2^n-1)$ by considering the two
intercalates:
$\{(0,2^n-1;2^n-1),(0,2^{n-1}-1;2^{n-1}-1),(2^{n-1},2^n-1,2^{n-1}-1),(2^{n-1},2^{n-1}-1;2^n-1)\}$
and
$\{(2^n-1,0;2^n-1),(2^n-1,2^{n-1},2^{n-1}-1),(2^{n-1}-1,0;2^{n-1}-1),(2^{n-1}-1,2^{n-1},2^{n-1},2^n-1)\}$.
For the entries $\{ (x,x;0) | 0 \leq x \leq 2^{n-2}-1 \}$, the
following intercalate proves the necessity of each entry:
$\{(x,x;0),(x,2^{n-2};x+2^{n-2}),(2^{n-2},x;x+2^{n-2}),(2^{n-2},2^{n-2};0)\}$.
For any other entry $x \in D$, the intercalate which suffices to
prove the necessity of $x$ in $C$ also suffices for $D$. Therefore
$D$ is a critical set of the required size. }\end{proof}

\section{The elementary abelian 3-group of order 9}
We create an integer program $IP_1$ similar to the above by adding
constraints based on all the Latin trades on 6 or less rows,
columns, or elements which have size less than or equal to 20. We
found, using WSAT, the following critical set $C_{29}$ of size 29
in the elementary abelian 3-group of order 9, $L(9)$.  Donovan,
Cooper, Nott and Seberry \cite{dcns} had found that any critical
set $L(9)$ contained at least 21 entries.

\begin{center}
$
\begin{latinsq}[1]
\hline &&&&&&&&8 \\
\hline &2&&4&&&7&&6 \\
\hline &&1&&3&&&&7 \\
\hline 3&4&&&&&0&& \\
\hline 4&5&&&&6&&2& \\
\hline &&&8&&&2&0& \\
\hline 6&&8&&1&&&&5 \\
\hline &&6&&&0&&5&3 \\
\hline &&7&2&&1&&& \\
\hline
\end{latinsq}
$
\end{center}

We create a new integer program $IP_2$ with all 324 Latin trades
of size 6 in the Latin square, and add symmetry considerations as
follows. We call rows 0 to 2, 3 to 5, and 6 to 8 of $L(9)$ $R_1$,
$R_2$ and $R_3$ respectively, and similarly call columns 0 to 2, 3
to 5, and 6 to 8 $C_1$, $C_2$ and $C_3$. $E_1$, $E_2$, and $E_3$
represent elements 0 to 2, 3 to 5, and 6 to 8. There are nine
subsquares defined by the row and column sets.  If $C$ is a
critical set in $L(9)$ we refer to these subsquares as
$C(R_i,C_j)$.

\begin{theorem}
If $C$ is a critical set in $L(9)$ then there exists another
critical set $D$ in $L(9)$ such that $|D(R_1,C_1)| \geq
|D(R_i,C_j)|$ for all $(i,j) \in \{1,2,3\}\times\{1,2,3\}$.
\end{theorem}

\begin{proof}{
In $C$, we can permute $R_1$ with $R_2$ or
$R_3$, and $C_1$ with $C_2$ or $C_3$ until the square defined by
$R_1$ and $C_1$ contains at least as many entries as any of the
other eight defined subsquares.  Call this property $RC$.

Each of these permutations is based on the fact that for any
critical set $C$ in $L(9)$ there exists an isomorphic critical set
in $L(9)$ after any of these permutations are performed.  This is
true because after the permutations above, we can permute $R_2$
with $R_3$, $C_2$ with $C_3$, or $E_1$, $E_2$ and $E_3$ to produce
a critical set isomorphic to $C$ in $L(9)$.  For example, for the
critical set $C_{29}$ above, we can swap $R_1$ with $R_3$, then we
can swap $E_1$ with $E_3$ and $C_2$ with $C_3$ to produce a
critical set in $L(9)$ with the property $RC$.  If for another
set, we swapped $C_1$ with $C_3$, we could then swap $R_2$ with
$R_3$ and $E_1$ with $E_3$ to obtain a critical set in $L(9)$.
Similar arguments apply for other permutations of rows and/or
columns.}\end{proof}

Thus $IP_2$ has another eight constraints for the property $RC$,
as in Theorem 2. Solving $IP_2$ with CPLEX, we find that any
critical set in this square is of size at least 24, improving the
bound of Donovan et al.

\begin{theorem}
If $C$ is a critical set of minimal size in $L(3^n)$, $n \geq 2$,
$|C| \geq 24.9^{n-2}$.
\end{theorem}

The inability of WSAT to find a solution of size less than 26
strongly suggests that 26 is the size of the smallest solution to
$IP_2$.

Solutions of size 28 to the integer program $IP_1$ above have been
found; some such solutions have trades with more than 6 rows,
columns, and elements missing. Thus a better approach, perhaps
based on one of the three suggestions above, is needed to prove
that 29 is the optimal solution, or that a better solution can be
found.  It would be symmetrically aesthetic if a critical set of
size 27 could be found with exactly 3 entries from each of the 36
$3 \times 3$ subsquares here. Over 540,000 solutions can be found
for an integer program with precisely these constraints
(constraints based on the 324 trades of size 6, plus constraints
ensuring there are exactly three entries in each subsquare).

\begin{conjecture}
The critical set of minimal size in $L(9)$ is of size 29.
\end{conjecture}

\end{document}